\documentclass [10pt] {article}

 \usepackage [cp1251] {inputenc}
 \usepackage [english] {babel}
 \usepackage {indentfirst}
 \usepackage {amssymb}
 \usepackage {stmaryrd}
 \usepackage [mathscr] {eucal}
 \usepackage {mathtools}
 \usepackage {amsthm}
 \usepackage [curve] {xypic}
 \usepackage {hyperref}

 \let \til\~           
 \let \gra\`
 \let \acu\'

 \def \?{{?}}
 \def \6{$'$}

 \def \Z{{\mathbb Z}}
 \def \Q{{\mathbb Q}}

 \def \E{{\boldsymbol{\mathcal E}}}            

 \def \F{{\mathscr F}} 

 \def \AA{{\boldsymbol A}}         
 \def \U{{\mathscr U}}             
 \def \FF{{\boldsymbol F}}         
 \def \N{{\mathscr N}}             

 \def \p{{\mathbf p}}   
 \def \r{{\mathbf R}}

 \def \id{{\mathrm{id}}}
 \def \pro{{\mathrm{pr}}}          
 \def \ins{{\mathrm{in}}}          
 \def \inc{{\mathrm{in}}}          
 \def \oBou{{({\Bou})}}            

 \def \O{{\displaystyle\raisebox{0.99ex}{$\scriptscriptstyle\boldsymbol<$}\mkern-4.1mu{|}}}
 \def \o{{\scriptstyle\raisebox{0.51ex}{$\scriptscriptstyle<$}\mkern-4.4mu{|}}}
 \def \0{{\mathchoice\O\O\o\o}}

 \def \<{\langle}      
 \def \>{\rangle}

 \def \LQ{$\scriptscriptstyle<$}                 
 \def \`{{\mathchoice
 {\raisebox{0.20ex}\LQ}
 {\raisebox{0.20ex}\LQ}
 {\raisebox{0.10ex}\LQ}
 {\raisebox{0.05ex}\LQ}
 }}

 \def \RQ{$\scriptscriptstyle>$}                 
 \def \'{{\mathchoice
 {\raisebox{0.20ex}\RQ}
 {\raisebox{0.20ex}\RQ}
 {\raisebox{0.10ex}\RQ}
 {\raisebox{0.05ex}\RQ}
 }}

 \def \({\llbracket}             
 \def \){\rrbracket}
 \def \[{[}            
 \def \]{]}

 \def \pp{\prescript+{}}           

 \DeclareMathOperator {\Ker} {Ker}
 \let \Im\undefined
 \DeclareMathOperator {\Im} {Im}
 \DeclareMathOperator {\ord} {ord}
 \let \deg\undefined
 \DeclareMathOperator {\deg} {deg}

 \def \bigbou{\bigvee}             

 \newcommand* {\bds} [3] {         
 \vcenter {
 \hrule height #1
 \vspace {#2}
 \hbox{$#3$}
 \vspace {#2}
 \hrule height #1
 }
 }

 \def \sbigBou{{\mathchoice        
 {\bds{0.20ex}{0.23ex}{\displaystyle\bigvee}}
 {\bds{0.15ex}{0.25ex}{\textstyle\bigvee}}
 {\bds{0.13ex}{0.21ex}{\scriptstyle\bigvee}}
 {\bds{0.09ex}{0.15ex}{\scriptscriptstyle\bigvee}}
 }}
 \def \bigBou{\mathop{\sbigBou}}

 \def \bou{\vee}       
 \def \Bou{\mathbin{\overline{\underline\vee}}}
 \def \cro{\times}     
 \def \Cro{\mathbin{\overline{\underline\times}}}
 \def \cat{\mathbin*}  

 \def \le{\leqslant}
 \def \ge{\geqslant}

 \def \-{\overline}    
 \def \~{\widetilde}
 \def \^{\widehat}

 \def \+{\boldsymbol}  

 \def \xto{\xrightarrow}

 \newcommand* {\bin} [2]
 {\Big({{#1}\atop{#2}}\Big)}

 \newcommand* {\tbin} [2]
 {\big({{#1}\atop{#2}}\big)}

 \renewcommand* {\%} [2]
 {\overset{#2}#1}

 \renewcommand* {\$} [3]
 {\overset{#2}{\underset{#3}#1}}

 \renewcommand* {\|} [2]
 {\mathbin{{#1}|_{#2}}}

 \renewcommand* {\1} [1] 
 {(#1)}

 \newdir {:=}{{}}      

 \newcommand* {\lelt} [2]          
 {\ar@{:=}@(dl,ul)+/l#1/;+/l#1/^{#2}}

 \newcommand* {\relt} [2]          
 {\ar@{:=}@(dr,ur)+/r#1/;+/r#1/_{#2}}

 \newcommand* {\suphead} [1]
 {\section* {\centerline{\normalsize\sc #1}}}

 \newcommand* {\head} [1]
 {\subsubsection* {\mathversion{bold}#1}}

 \newenvironment* {claim} [1] []
 {\begin{trivlist}\item [\hskip\labelsep {\bf #1}] \it}
 {\end{trivlist} }

 \newenvironment* {demo} [1] []
 {\begin{trivlist}\item [\hskip\labelsep {\it #1}] }
 {\end{trivlist} }

 \makeatletter
 \def \QED{\displaymath@qed}
 \makeatother

\begin {document}

 \title {\large\bf
         Homotopy similarity of maps.
         Maps of the circle}

 \author {\normalsize\rm
          S.~S.~Podkorytov}

 \date {}

 \maketitle

 \vspace {-2\bigskipamount}
 
 \begin {abstract} \noindent
 We describe
 the relation of $r$-similarity and
 finite-order invariants
 on the homotopy set $[S^1,Y]=\pi_1(Y)$.
 \end {abstract}


 \head {\S~1. Introduction}


 This paper continues \cite{sim}.
 We adopt notation and conventions thereof.
 Here we are mainly interested in the set $[S^1,Y]=\pi_1(Y)$;
 in Part~I,
 however,
 we consider a more general case.
 Let $X$ and $Y$ be cellular spaces,
 $X$ compact.
 Let $X$ be equipped with maps
 $\mu:X\to X\bou X$ (comultiplication) and
 $\nu:X\to X$ (coinversion).
 The set $Y^X$ carries the operations
 $$
 (a,b)
 \mapsto
 (a\cat b:X\xto{\mu}X\bou X\xto{a\Bou b}Y)
 $$
 and
 $$
 a
 \mapsto
 (a^\dag:X\xto{\nu}X\xto{a}Y).
 $$
 We suppose that
 the set $[X,Y]$ is a group with
 the identity $1=[\0^X_Y]$,
 the multiplication
 $$
 [a][b]=[a\cat b],
 $$
 and
 the inversion
 $$
 [a]^{-1}=[a^\dag].
 $$
 Under these assumptions,
 we call $(X,\mu,\nu;Y)$ an {\it admissible couple}.

 Put
 $$
 [X,Y]^{(r+1)}=
 \{\,\+a\in[X,Y]\mid1\%\sim r\+a\,\}.
 $$
 We get the filtration
 $$
 [X,Y]=
 [X,Y]^{(1)}\supseteq
 [X,Y]^{(2)}\supseteq\dotso.
 $$
 We prove that
 the subsets $[X,Y]^{(r+1)}$
 are normal subgroups and
 form an N-series
 (Theorems 4.1 and 4.3).
 The equivalence holds
 $$
 \+a\%\sim r\+b
 \quad
 \Leftrightarrow
 \quad
 \+a^{-1}\+b\in[X,Y]^{(r+1)}
 $$
 (Theorem 4.2).

 In Part~III,
 we concentrate on the case $X=S^1$
 (with the standard $\mu$ and $\nu$),
 when $[X,Y]=\pi_1(Y)$.
 We prove that
 $$
 \pi_1(Y)^{(r+1)}=\gamma^{r+1}\pi_1(Y)
 $$
 (Theorem~11.2).
 Here,
 as usual,
 $$
 G=
 \gamma^{1}G\supseteq
 \gamma^{2}G\supseteq\dotso
 $$
 is the lower central series of a group $G$.

 For a homotopy invariant
 (i.~e., a function)
 $f:\pi_1(Y)\to L$,
 where $L$ is an abelian group,
 its order $\ord f\in\{-\infty,0,1,\dotsc,\infty\}$
 is defined
 (see \cite[\S~1]{sim}).
 We prove that
 $\ord f=\deg f$
 (Theorem~12.2).
 Recall that,
 for a function $f:G\to L$,
 where $G$ is a group,
 its degree $\deg f$ is defined
 (see \S~12).

 Do invariants of order at most $r$ distinguish
 elements of $\pi_1(Y)$
 that are not $r$-similar?
 In general,
 no.
 For $r\ge3$,
 there is
 a group $G$ and
 an element $g\in G\setminus\gamma^{r+1}G$
 such that,
 for any
 abelian group $L$ and
 function $f:G\to L$ of degree at most $r$,
 one has $f(1)=f(g)$
 (see
 \cite{R} for $r=3$
 and
 \cite[Ch.~2]{MP}).
 Take a cellular space $Y$ with $\pi_1(Y)=G$.
 Then,
 by Theorems 11.2 and 12.2,
 the homotopy classes $1$ and $g$ in $\pi_1(Y)$
 are not $r$-similar,
 but cannot be distinguished by invariants
 of order at most $r$.

 In Part~II,
 which does not depend on the rest of the paper,
 we prove group-theoretic Theorem~9.1,
 which we need
 for the proof of the above-mentioned Theorem~11.2.


 \suphead {Part I}


 In this part,
 we
 discuss operations over coherent ensembles of maps between
 arbitrary spaces
 (\S\S\ 2 and 3) and
 give our results concerning an arbitrary admissible couple
 (\S~4).


 \head {\S~2. Compositions}

 Let
 $X$, $Y$, $X'$, and $Y'$ be spaces and
 $k:X'\to X$ and $h:Y\to Y'$ be maps.
 Introduce the homomorphisms
 $$
 k^\#:
 \<Y^X\>\to
 \<Y^{X'}\>,
 \qquad
 \`a\'\mapsto\`a\circ k\',
 $$
 and
 $$
 h_\#:
 \<Y^X\>\to
 \<{Y'}^X\>,
 \qquad
 \`a\'\mapsto\`h\circ a\'.
 $$

 \begin {claim} [2.1. Lemma.]
 We have
 $$
 k^\#(\<Y^X\>^{(r+1)})\subseteq
 \<Y^{X'}\>^{(r+1)}
 \quad
 \text{and}
 \quad
 h_\#(\<Y^X\>^{(r+1)})\subseteq
 \<{Y'}^X\>^{(r+1)}.
 $$
 \end {claim}

 \begin {demo} [Proof.]
 Take an ensemble $A\in\<Y^X\>^{(r+1)}$.

 To show that $k^\#(A)\in\<Y^{X'}\>^{(r+1)}$,
 we
 take $T'\in\F_{r}(X')$ and
 check that $k^\#(A)|_{T'}=0$.
 We have the commutative diagram
 $$
 \xymatrix {
 \<Y^{X'}\>
 \lelt{2ex}{k^\#(A)}
 \ar[d]_-{\?|_{T'}}
 &&
 \<Y^X\>
 \relt{2ex}{A}
 \ar[ll]_-{k^\#}
 \ar[d]^-{\?|_{k(T')}}
 \\
 \<Y^{T'}\>
 \lelt{2ex}{k^\#(A)|_{T'}}
 &&
 \<Y^{k(T')}\>
 ,
 \relt{2ex}{A|_{k(T')}=0}
 \ar[ll]_-{q^\#}
 }
 $$
 where $q=k|_{T'\to k(T')}$.
 Since $k(T')\in\F_r(X)$,
 we have $A|_{k(T')}=0$.
 By the diagram,
 $k^\#(A)|_{T'}=0$.

 To show that $h_\#(A)\in\<{Y'}^X\>^{(r+1)}$,
 we
 take $T\in\F_{r}(X)$ and
 check that $h_\#(A)|_T=0$.
 We have the commutative diagram
 $$
 \xymatrix {
 \<Y^X\>
 \lelt{2ex}{A}
 \ar[rr]^-{h_\#}
 \ar[d]_-{\?|_T}
 &&
 \<{Y'}^X\>
 \relt{2ex}{h_\#(A)}
 \ar[d]^-{\?|_T}
 \\
 \<Y^T\>
 \lelt{2ex}{0=A|_T}
 \ar[rr]^-{h_\#}
 &&
 \<{Y'}^T\>
 .
 \relt{2ex}{h_\#(A)|_T}
 }
 $$
 We have $A|_T=0$.
 By the diagram,
 $h_\#(A)|_T=0$.
 \qed
 \end {demo}

 \begin {claim} [2.2. Corollary.]
 Let $a,b\in Y^X$ satisfy $a\%\sim rb$.
 Then
 $a\circ k\%\sim rb\circ k$ in ${X'}^Y$ and
 $h\circ a\%\sim rh\circ b$ in $X^{Y'}$.
 \end {claim}

 \begin {demo} [Proof.]
 There is an ensemble $A\in\<Y^X\>$,
 $$
 A=
 \sum_i
 u_i\`a_i\',
 $$
 where $a_i\sim a$,
 such that
 $A\%=r\`b\'$.
 By Lemma~2.1,
 $k^\#(A)\%=r\`b\circ k\'$ and
 $h_\#(A)\%=r\`h\circ b\'$.
 Since all the maps of $k^\#(A)$ are homotopic to $a\circ k$,
 we get $a\circ k\%\sim rb\circ k$.
 Since all the maps of $h_\#(A)$ are homotopic to $h\circ b$,
 we get $h\circ a\%\sim rh\circ b$.
 \qed
 \end {demo}


 \head {\S~3. Joining coherent ensembles}

 Let $X_1$, $X_2$, and $Y$ be spaces.
 Introduce the homomorphism
 $$
 \oBou:
 \<Y^{X_1}\>\otimes\<Y^{X_2}\>
 \to
 \<Y^{X_1\bou X_2}\>,
 \qquad
 \`a\'\otimes\`b\'
 \mapsto
 \`a\Bou b\'.
 $$

 \begin {claim} [3.1. Lemma.]
 For $p,q\ge0$,
 we have
 $$
 \oBou(\<Y^{X_1}\>^{(p)}\otimes\<Y^{X_2}\>^{(q)})
 \subseteq
 \<Y^{X_1\bou X_2}\>^{(p+q)}.
 $$
 \end {claim}

 \begin {demo} [Proof.]
 Take
 $A\in\<Y^{X_1}\>^{(p)}$ and
 $B\in\<Y^{X_2}\>^{(q)}$.
 We show that
 $\oBou(A\otimes B)\in\<Y^{X_1\bou X_2}\>^{(p+q)}$.
 Take $T\in\F_{p+q-1}(X_1\bou X_2)$.
 We check that
 $\oBou(A\otimes B)|_T=0$.
 We have $T=T_1\bou T_2$
 for some finite subspaces $T_i\subseteq X_i$,
 $i=1,2$.
 We have the commutative diagram
 $$
 \xymatrix {
 \<Y^{X_1}\>\otimes\<Y^{X_2}\>
 \lelt{6ex}{A\otimes B}
 \ar[rr]^-{\oBou}
 \ar[d]_-{\?|_{T_1}\otimes\?|_{T_1}}
 &&
 \<Y^{X_1\bou X_2}\>
 \relt{4ex}{\oBou(A\otimes B)}
 \ar[d]^-{\?|_T}
 \\
 \<Y^{T_1}\>\otimes\<Y^{T_2}\>
 \lelt{6ex}{A|_{T_1}\otimes B|_{T_2}}
 \ar[rr]^-{\oBou}
 &&
 \<Y^T\>
 .
 \relt{4ex}{\oBou(A\otimes B)|_T}
 }
 $$
 We have
 $T_1\in\F_{p-1}(X_1)$ or
 $T_2\in\F_{q-1}(X_2)$.
 Thus
 $A|_{T_1}=0$ or
 $B|_{T_2}=0$.
 By the diagram,
 $\oBou(A\otimes B)|_T=0$.
 \qed
 \end {demo}


 \head {\S~4. Similarity for an admissible couple}

 Let $(X,\mu,\nu;Y)$ be an admissible couple.

 \begin {claim} [4.1. Theorem.]
 $[X,Y]^{(r+1)}\subseteq[X,Y]$ is a subgroup.
 \end {claim}

 \begin {demo} [Proof.]
 To show that
 $[X,Y]^{(r+1)}$ is closed under multiplication,
 we take $a,b\in Y^X$ such that
 $\0\%\sim ra$ and
 $\0\%\sim rb$
 and
 check that
 $\0\%\sim ra\cat b$.
 There are ensembles $D,E\in\<Y^X\>$,
 $$
 D=
 \sum_i
 u_i\`d_i\'
 \quad
 \text{and}
 \quad
 E=
 \sum_j
 v_j\`e_j\',
 $$
 where
 $d_i\sim\0$ and
 $e_j\sim\0$,
 such that
 $D\%=r\`a\'$ and
 $E\%=r\`b\'$.
 Consider
 the maps $a\Bou b,d_i\Bou e_j:X\bou X\to Y$ and
 the ensemble $F\in\<Y^{X\bou X}\>$,
 $$
 F=
 \sum_{i,j}
 u_iv_j\`d_i\Bou e_j\'.
 $$
 We have
 \begin {multline*}
 \`a\Bou b\'-F=
 \oBou(\`a\'\otimes\`b\')-\oBou(D\otimes E)=
 \\ =
 \oBou((\`a\'-D)\otimes\`b\')+\oBou(D\otimes(\`b\'-E)
 \in\<Y^{X\bou X}\>^{(r+1)},
 \end {multline*}
 where $\in$ holds
 by Lemma~3.1.
 Since all the maps of $F$ are null-homotopic,
 we get $\0\%\sim ra\Bou b$.
 Since $a\cat b=(a\Bou b)\circ\mu$,
 Corollary~2.2 yields $\0\%\sim ra\cat b$.

 Take $a\in Y^X$ such that
 $\0\%\sim ra$.
 Since $a^\dag=a\circ\nu$,
 Corollary~2.2 yields $\0\%\sim ra^\dag$.
 Thus
 $[X,Y]^{(r+1)}$ is closed under inversion.
 \qed
 \end {demo}

 \begin {claim} [4.2. Theorem.]
 For $\+a,\+b\in[X,Y]$,
 we have
 \begin {equation} \label {cri}
 \+a\%\sim r\+b
 \quad
 \Leftrightarrow
 \quad
 \+a^{-1}\+b\in[X,Y]^{(r+1)}.
 \end {equation}
 \end {claim}

 \begin {demo} [Proof.]
 It suffices to check the implication
 $$
 a\%\sim rb
 \quad
 \Rightarrow
 \quad
 c\cat a\%\sim rc\cat b
 $$
 for $a,b,c\in Y^X$.
 Given an ensemble $A\in\<Y^X\>$,
 $$
 A=
 \sum_i
 u_i\`a_i\',
 $$
 where $a_i\sim a$,
 such that
 $A\%=r\`b\'$,
 consider the ensemble $F\in\<Y^{X\bou X}\>$,
 $$
 F=
 \sum_i
 u_i\`c\Bou a_i\'.
 $$
 We have
 $$
 \`c\Bou b\'-F=
 \oBou(\`c\'\otimes(\`b\'-A))
 \in\<Y^{X\bou X}\>^{(r+1)},
 $$
 where $\in$ holds by Lemma~3.1.
 Thus
 $F\%=r\`c\Bou b\'$.
 Since $c\Bou a_i\sim c\Bou a$,
 we get $c\Bou a\%\sim rc\Bou b$.
 Taking composition with $\mu$,
 we get $c\cat a\%\sim rc\cat b$
 by Corollary~2.2.
 \qed
 \end {demo}

 Theorems 4.1 and 4.2 imply that
 the relaton $\%\sim r$ on $[X,Y]$ is an equivalence,
 which is a special case of \cite[Theorem~8.1]{sim}
 (note that
 we did not use it here).

 One can prove similarly that
 \begin {equation} \label {cri-1}
 \+a\%\sim r\+b
 \quad
 \Leftrightarrow
 \quad
 \+b\+a^{-1}\in[X,Y]^{(r+1)}.
 \end {equation}
 It follows from \eqref{cri} and \eqref{cri-1} that
 the subgroup $[X,Y]^{(r+1)}\subseteq[X,Y]$ is normal.
 This is a special case of the following theorem.

 Let $\(\ ,\ \)$ denote the group commutator.

 \begin {claim} [4.3. Theorem.]
 Put $M^s=[X,Y]^{(s)}\subseteq[X,Y]$.
 Then $\(M^p,M^q\)\subseteq M^{p+q}$.
 \end {claim}

 \begin {demo} [Proof.]
 Introduce the map
 $$
 \zeta:
 X
 \xto{\mu^{(3)}}
 X\bou X\bou X\bou X
 \xto{(\ins_1\circ\nu)\Bou(\ins_2\circ\nu)\Bou\ins_1\Bou\ins_2}
 X\bou X,
 $$
 where
 $$
 \mu^{(3)}:
 X
 \xto{\mu}
 X\bou X
 \xto{\mu\bou\id_X}
 X\bou X\bou X
 \xto{\mu\bou\id_X\bou\id_X}
 X\bou X\bou X\bou X
 $$
 ($4$-fold comultiplication).
 For $a,b\in Y^X$,
 we have
 \begin {equation} \label {ctr}
 [(a\Bou b)\circ\zeta]=\([a],[b]\)
 \end {equation}
 in the group $[X,Y]$.

 Take $a,b\in Y^X$ such that
 $\0\%\sim{p-1}a$ and
 $\0\%\sim{q-1}b$.
 We show that
 $\0\%\sim{p+q-1}(a\Bou b)\circ\zeta$.

 There are ensembles $D,E\in\<Y^X\>$,
 $$
 D=
 \sum_i
 u_i\`d_i\'
 \quad
 \text{and}
 \quad
 E=
 \sum_j
 v_j\`e_j\',
 $$
 where
 $d_i\sim\0$ and
 $e_j\sim\0$,
 such that
 $D\%={p-1}\`a\'$ and
 $E\%={q-1}\`b\'$.
 Consider
 the ensemble $F\in\<Y^{X\bou X}\>$,
 $$
 F=
 \sum_i
 u_i\`d_i\Bou b\'+
 \sum_j
 v_j\`a\Bou e_j\'-
 \sum_{i,j}
 u_iv_j\`d_i\Bou e_j\'.
 $$
 We have
 $$
 \`a\Bou b\'-F=
 \oBou((\`a\'-D)\otimes(\`b\'-E))
 \in\<Y^{X\bou X}\>^{(p+q)},
 $$
 where $\in$ holds by Lemma~3.1.
 Thus
 $F\%={p+q-1}\`a\Bou b\'$.
 By Lemma~2.1,
 $\zeta^\#(F)\%={p+q-1}\`(a\Bou b)\circ\zeta\'$.
 By \eqref{ctr},
 all the maps of $\zeta^\#(F)$ are null-homotopic.
 Thus
 we get $\0\%\sim{p+q-1}(a\Bou b)\circ\zeta$.
 \qed
 \end {demo}


 \suphead {Part II}


 In this part,
 which
 is algebraic and
 does not depend on the rest of the paper,
 we prove Theorem~9.1.


 \head {\S~5. Cultured sets}

 Let $E$ be a set.
 Consider the $\Q$-algebra $\Q^E$ of functions $E\to\Q$.
 A {\it culture\/} on $E$ is a filtration $\Phi=(\Phi_s)_{s\ge0}$
 of $\Q^E$
 by $\Q$-submodules
 $$
 \Phi_0\subseteq
 \Phi_1\subseteq
 \dotso\subseteq
 \Q^E
 $$
 such that
 $$
 1\in\Phi_0
 \quad
 \text{and}
 \quad
 \Phi_s\Phi_t\subseteq
 \Phi_{s+t}.
 $$
 A set equipped with a culture is called a {\it cultured set}.
 The culture of a cultured set $E$ is denoted by $\Phi^E$.

 A way to define a culture on a set $E$ is
 to choose a collection of pairs $(u_i,s_i)$,
 where
 $u_i\in\Q^E$ is a function and
 $s_i\ge1$ is a number called the {\it weight}, and
 to let $\Phi_s$ be spanned by all products $u_{i_1}\dotso u_{i_p}$
 ($p\ge0$)
 with $s_{i_1}+\dotso+s_{i_p}\le s$.
 We define the cultured set
 \begin {equation} \label {Qms}
 \Q^m_{s_1\dotso s_m}
 \end {equation}
 as $\Q^m$ with the culture given the collection $(\xi_i,s_i)$,
 $i\in\1m$,
 where $\xi_i:\Q^m\to\Q$ is the $i$th coordinate.
 Hereafter,
 we put $\1m=\{1,\dotsc,m\}$.
 The cultured set
 $$
 \Z^m_{s_1\dotso s_m}
 $$
 is defined similarly.
 We put
 $\Q_s=\Q^1_s$ and
 $\Z_s=\Z^1_s$.

 A function $g:E\to F$ between cultured sets
 is called a {\it cultural morphism\/} if
 the induced algebra homomorphism $h^\#:\Q^F\to\Q^E$
 satisfies $g^\#(\Phi^F_s)\subseteq\Phi^E_s$
 for all $s$.
 A function
 $$
 g:
 \Q^m_{s_1\dotso s_m}\to
 \Q^n_{t_1\dotso t_n}
 $$
 is a cultural morphism
 if and only if
 it has the form
 $$
 g(x_1,\dotsc,x_m)=
 (P_j(x_1,\dotsc,x_m))_{j\in\1n},
 $$
 where
 $P_j$ is a rational polynomial of degree at most $t_j$
 with respect to its arguments having weights $s_1,\dots,s_m$.
 Cultural maps
 $$
 \Z^m_{s_1\dotso s_m}\to
 \Z^n_{t_1\dotso t_n}
 $$
 are characterized similarly
 (their coordinate polynomials need not have integer coefficients).

 Cultured sets and cultural morphisms form a category with products.
 We have
 $$
 \Q^m_{s_1\dotso s_m}\cro
 \Q^n_{t_1\dotso t_n}=
 \Q^{m+n}_{s_1\dotso s_mt_1\dotso t_n}
 \quad
 \text{and}
 \quad
 \Z^m_{s_1\dotso s_m}\cro
 \Z^n_{t_1\dotso t_n}=
 \Z^{m+n}_{s_1\dotso s_mt_1\dotso t_n}.
 $$

 A cultural morphism $g:E\to F$ is called a {\it cultural immersion\/} if
 $g^\#(\Phi^F_s)=\Phi^E_s$
 for all $s$.
 Then
 a function $f:D\to E$,
 where $D$ is a cultured set,
 is a cultural morphism if
 the composition
 $$
 D
 \xto{f}
 E
 \xto g
 F
 $$
 is.
 If the composition
 $$
 E
 \xto g
 F
 \xto h
 G
 $$
 of two cultural morphisms is a cultural immersion,
 then $g$ is.


 \head {\S~6. The truncated free algebra $\AA/\AA^{(r+1)}$}

 Consider the algebra $\AA$ of rational polynomials
 in non-commuting variables $T_1,\dotsc,T_n$.
 It is graded in the standard way,
 $$
 \AA=
 \bigoplus_{s\ge0}
 \AA_s.
 $$
 Introduce the ideals $\AA^{(s)}\subseteq\AA$,
 $$
 \AA^{(s)}=
 \bigoplus_{t\ge s}
 \AA_t.
 $$

 We
 fix $r\ge0$ and
 consider
 the algebra $\AA/\AA^{(r+1)}$.
 Let $\-T_i\in\AA/\AA^{(r+1)}$ be the image of $T_i$.
 An element $w\in\AA/\AA^{(r+1)}$ has the form
 $$
 w=
 \sum_{s\ge0,\,i_1,\dotsc,i_s}
 w^{(s)}_{i_1\dotsc i_s}\-T_{i_1}\dotsc\-T_{i_s},
 $$
 where $w^{(s)}_{i_1\dotsc i_s}\in\Q$ are
 uniquely defined
 for $s\le r$
 and
 arbitrary
 for greater $s$.
 Introduce
 the group
 $$
 \U=
 1+\AA^{(1)}/\AA^{(r+1)}\subseteq
 (\AA/\AA^{(r+1)})^\times
 $$
 with the filtration by subgroups
 $$
 \U=
 \U^{(1)}\supseteq
 \U^{(2)}\supseteq\dotso
 $$
 where
 $$
 \U^{(s)}=
 1+\AA^{(s)}/\AA^{(r+1)},
 \qquad
 s\le r+1,
 $$
 and
 $\U^{(s)}=1$
 for $s\ge r+1$.
 We have $\(\U^{(s)},\U^{(t)}\)\subseteq\U^{(s+t)}$.
 In particular,
 $\gamma^s\U\subseteq\U^{(s)}$.

 We equip the set $\U$ with the culture
 given by the collection of pairs
 $$
 (\xi^{(s)}_{i_1\dotso i_s},s),
 \qquad
 1\le s\le r,
 \quad
 i_1,\dotsc,i_s\in\1n,
 $$
 where
 $$
 \xi^{(s)}_{i_1\dotso i_s}:\U\to\Q,
 \qquad
 u\mapsto u^{(s)}_{i_1\dotso i_s}:\U\to\Q
 $$
 for
 $$
 u=
 \sum_{t\ge0,\,j_1,\dotsc,j_t}
 u^{(t)}_{j_1\dotsc j_t}\-T_{j_1}\dotsc\-T_{j_t}
 $$
 with $u^{(0)}=1$.
 Clearly,
 the cultured set $\U$ is a special case of \eqref{Qms}.

 For a group $G$,
 let $M_G:G\cro G\to G$ be the multiplication.

 \begin {claim} [6.1. Lemma.]
 The function
 $$
 M_\U:
 \U\cro\U\to
 \U
 $$
 is a cultural morphism.
 \end {claim}

 \begin {demo} [Proof.]
 Given
 $u,v\in\U$,
 $$
 u=
 \sum_{s\ge0,\,i_1,\dotsc,i_s}
 u^{(s)}_{i_1\dotsc i_s}\-T_{i_1}\dotsc\-T_{i_s}
 \quad
 \text{and}
 \quad
 v=
 \sum_{t\ge0,\,j_1,\dotsc,j_s}
 v^{(t)}_{j_1\dotsc j_t}\-T_{j_1}\dotsc\-T_{j_t}
 $$
 with $u^{(0)}=v^{(0)}=1$,
 we have
 $$
 uv=
 \sum_{\substack{
         s\ge0,\,i_1,\dotsc,i_s, \\
         t\ge0,\,j_1,\dotsc,j_t
 }}
 u^{(s)}_{i_1\dotsc i_s}
 v^{(t)}_{j_1\dotsc j_t}
 \-T_{i_1}\dotsc\-T_{i_s}
 \-T_{j_1}\dotsc\-T_{j_t}.
 $$
 We consider
 $u^{(s)}_{i_1\dotsc i_s}$
 and
 $v^{(t)}_{j_1\dotsc j_t}$
 with $s,t>0$
 here as variables of weights
 $s$
 and
 $t$,
 respectively.
 In the last expression,
 the monomial in $\-T_i$ has degree $s+t$,
 and
 its coefficient
 is $u^{(s)}_{i_1\dotsc i_s}v^{(t)}_{j_1\dotsc j_t}$
 and thus
 has degree $s+t$.
 Thus
 the total coefficient
 of each monomial in $\-T_i$ of some degree $z>0$
 in this series
 is a polynomial in
 $u^{(s)}_{i_1\dotsc i_s}$
 and
 $v^{(t)}_{j_1\dotsc j_t}$
 of degree at most $z$.
 \qed
 \end {demo}

 \begin {claim} [6.2. Lemma.]
 For $u\in\U^{(s)}$,
 the function
 $$
 E_u:
 \Z_s\to
 \U,
 \qquad
 x\to u^x,
 $$
 is a cultural morphism.
 Moreover,
 it extends to a cultural morphism
 $$
 \^E_u:
 \Q_s\to
 \U.
 $$
 \end {claim}

 \begin {demo} [Proof.]
 We have $u=1+w$
 in $\AA/\AA^{(r+1)}$
 for some $w\in\AA^{(s)}/\AA^{(r+1)}$,
 $$
 w=
 \sum_{t\ge s,\,i_1,\dotsc,i_t}
 w^{(t)}_{i_1\dotsc i_t}\-T_{i_1}\dotsc\-T_{i_t}.
 $$
 For $x\in\Z$,
 we have
 \begin {multline*}
 E_u(x)=
 u^x=
 (1+w)^x=
 \sum_{p\ge0}
 \bin xpw^p=
 \\
 =
 \sum_{p\ge0}\
 \sum_{\substack{
         t_1\ge s,\,i_{11},\dotsc,i_{1t_1},\\
         \dotsc,\\
         t_p\ge s,\,i_{p1},\dotsc,i_{pt_p}
 }}
 w^{(t_1)}_{i_{11},\dotsc,i_{1t_1}}
 \dotso
 w^{(t_p)}_{i_{p1},\dotsc,i_{pt_p}}
 \bin xp
 \cdot
 \\
 \cdot
 \-T_{i_{11}}\dotso\-T_{i_{1t_1}}
 \dotso
 \-T_{i_{p1}}\dotso\-T_{i_{pt_p}}.
 \end {multline*}
 Consider $x$ here as a variable of weight $s$.
 In the last expression,
 the monomial in $\-T_i$ has degree $t_1+\dotso+t_p$,
 which is at least $ps$.
 Its coefficient is
 a rational multiple of $\tbin xp$
 and thus
 a polynomial in $x$ of degree at most $ps$.
 Thus
 the total coefficient
 of each monomial in $\-T_i$ of some degree $z>0$
 in this series
 is a polynomial in $x$ of degree at most $z$.

 The extension $\^E_u$ exists automatically.
 \qed
 \end {demo}


 \head {\S~7. The free nilpotent group $\N$}

 Recall that
 we fix numbers
 $n$ and
 $r$.
 Let $\FF$ be the free group on the generators $Z_1,\dotsc,Z_n$.
 Consider the free nilpotent group $\N=\FF/\gamma^{r+1}\FF$.
 Put $\N^{(s)}=\gamma^s\N\subseteq\N$, $s\ge1$.

 Following Magnus,
 consider the homomorphism
 $$
 \rho:\N\to\U,
 \qquad
 \-Z_i\mapsto1+\-T_i.
 $$
 Hereafter,
 the bar denotes the projection to the proper quotient group.
 The homomorphism $\rho$ exists because
 $\gamma^{r+1}\U=1$,
 The quotient $\N^{(s)}/\N^{(s+1)}$ is
 abelian
 and
 finitely generated.
 Since
 $
 \rho(\N^{(s)})
 \subseteq
 \gamma^s\U
 \subseteq
 \U^{(s)}
 $,
 there is a homomorphism
 $$
 \sigma^{(s)}:
 \N^{(s)}/\N^{(s+1)}\to
 \AA_s
 $$
 such that
 $$
 \rho(h)\equiv
 1+\sigma^{(s)}(\-h)
 \pmod{\AA^{(s+1)}},
 \qquad
 h\in\N^{(s)}.
 $$
 By Magnus \cite{M}
 (see also \cite[Part~I, Ch.~IV, Theorem~6.3]{S}),
 $\N^{(s)}=\rho^{-1}(\U^{(s)})$.
 It follows that
 $\sigma^{(s)}$ are injective
 and
 $\N^{(s)}/\N^{(s+1)}$ are
 torsion-free
 and
 thus
 free abelian.
 It follows that
 there is a filtration
 $$
 \N=
 \N^1\supseteq
 \N^2\supseteq
 \dotso\supseteq
 \N^q\supseteq
 \N^{q+1}=1
 $$
 such that
 $\N^{(s)}=\N^{j_s}$
 for some $1=j_1\le\dotso\le j_{r+1}=q+1$
 and
 $\N^j/\N^{j+1}$ are infinite cyclic.
 For $j\le n+1$,
 we choose $\N^j$ be the subgroup generated by
 $\-Z_j,\dotsc,\-Z_n$
 and
 $\N^{(2)}$.
 Put $s_j=\max\{\,s\mid\,j_s\le j\,\}$,
 $j\in\1q$.
 Clearly,
 $1\le s_1\le\dotso\le s_q\le r$,
 $s_1=\dotso=s_n=1$,
 and
 $\N^j\subseteq\N^{(s_j)}$.
 The subgroups $\N^j\subseteq\N$ are normal.

 For each $j\in\1q$,
 choose an element $b_j\in\N^j$ such that
 $\-b_j$ generates $\N^j/\N^{j+1}$.
 In doing so,
 we put
 $$
 b_j=\-Z_j,
 \qquad
 j\in\1n.
 $$
 The collection $(b_1,\dotsc,b_q)$ is a ``Mal\6cev basis''
 \cite[4.2.2]{CMZ}.
 For $j\in\1{q+1}$,
 the function
 $$
 \beta^j:
 \Z^{q-j+1}
 \to
 \N^j,
 \qquad
 (x_j,\dotsc,x_q)
 \mapsto
 b_j^{x_j}\dotso b_q^{x_q},
 $$
 is bijective.
 We put
 $$
 \beta=
 \beta^1:
 \Z^q
 \to
 \N.
 $$
 The elements $\sigma^{(s_j)}(\-b_j)\in\AA$ are linearly independent.

 Any group $G$ carries the {\it immanent\/} culture $\Phi$
 with $\Phi_s$ consisting of all functions $G\to\Q$ of degree at most $s$
 (see \S~12).
 If $\N$ is equipped with its immanent culture,
 $$
 \beta:
 \Z^q_{s_1\dotso s_q}
 \to
 \N
 $$
 becomes a culture isomorphism.
 The proof is omitted.

 \begin {claim} [7.1. Lemma.]
 The composition
 $$
 \eta^j:
 \Z^{{q-j+1}}_{s_j\dotso s_q}
 \xto{\beta^j}
 \N^j
 \xto{\rho^j}
 \U,
 $$
 where $\rho^j=\rho|_{\N^j}$,
 is a cultural immersion.
 \end {claim}

 Introduce the projections
 $$
 \p:
 \Q^{q-j+1}
 \to
 \Q,
 \qquad
 (x_j\dotsc,x_q)\mapsto x_j,
 $$
 and
 $$
 \r:
 \Q^{q-j+1}
 \to
 \Q^{q-j},
 \qquad
 (x_j\dotsc,x_q)\mapsto(x_{j+1},\dotsc,x_q).
 $$

 \begin {demo} [Proof.]
 We show that
 $\eta^j$ is a cultural morphism
 by backward induction on $j$.
 For $j=q+1$,
 the assertion is trivial.
 Take $j\le q$.
 Since $b^j\in\N^{(s_j)}$,
 we have $\rho(b^j)\in\U^{(s_j)}$.
 We have the decomposition
 $$
 \eta^j:
 \Z^{q-j+1}_{s_j\dotso s_q}=
 \Z_{s_j}\cro\Z^{q-j}_{s_{j+1}\dotso s_q}
 \xto{E_{\rho(b_j)}\cro\eta^{j+1}}
 \U\cro\U
 \xto{M_\U}
 \U,
 $$
 where $E_{\rho(b_j)}:x\mapsto\rho(b_j)^x$.
 Here
 cultural morphisms are
 $E_{\rho(b_j)}$
 by Lemma~6.2,
 $\eta^{j+1}$
 by the induction hypothesis,
 and
 $M_\U$ by Lemma~6.1.
 Thus
 $\eta^j$ is a cultural morphism.

 For each $j\in\1q$,
 choose a linear functional $\phi_j:\AA_{s_j}\to\Q$ such that
 $\phi_j(\sigma^{(s_j)}(b_k))$ equals
 $1$
 for $k=j$ and
 $0$
 for all other $k$
 with $s_k=s_j$.
 Given
 $j\le q+1$
 and
 $x=(x_j,\dotsc,x_q)\in\Z^{q-j+1}$,
 we have
 $$
 \eta^j(x)=
 \rho(\beta^j(x))=
 \rho(b_j^{x_j}\dotso b_q^{x_q})=
 \rho(b_j^{x_j})\dotso\rho(b_q^{x_q}).
 $$
 Assume $j\le q$.
 Then
 $$
 \eta^j(x)\equiv
 1+
 \sum_{k\ge j:s_k=s_j}
 x_k\sigma^{(s_j)}(b_k)
 \pmod{\AA^{(s_j+1)}/\AA^{(r+1)}}
 $$
 in $\AA/\AA^{(r+1)}$
 and
 $$
 \eta^j(x)=
 \rho(b_j)^{x_j}\eta^{j+1}(\r(x)).
 $$
 Note that,
 for any linear functional $F:\AA_{s_j}\to\Q$,
 the composition
 $$
 F!:
 \U
 \xto\inc
 \AA/\AA^{(r+1)}
 \xto\pro
 \AA_{s_j}
 \xto F
 \Q_{s_j}
 $$
 is a cultural morphism.
 For $c\in\Q$,
 we have
 $$
 (c\phi_j)!\,(\eta^j(x))=cx_j,
 \qquad
 x=(x_j,\dotsc,x_q)\in\Z^{q-j+1}.
 $$
 
 We show that
 $\eta^j$ is a cultural immersion
 by constructing a cultural morphism
 $$
 \theta^j:
 \U\to
 \Q^{q-j+1}_{s_j\dotso s_q}
 $$
 such that
 $\theta^j\circ\eta^j$ is the inclusion
 $$
 \Z^{q-j+1}_{s_j\dotso s_q}\to
 \Q^{q-j+1}_{s_j\dotso s_q}.
 $$
 Backward induction on $j$.
 Let $\theta^{q+1}$ be the unique function $\U\to\Q^0$.
 Take $j\le q$.
 Introduce the cultural morphism
 $$
 L:
 \U
 \xto{(-\phi_j)!}
 \Q_{s_j}
 \xto{\^E_{\rho(b_j)}}
 \U
 $$
 where $\^E_{\rho(b_j)}$ is given by Lemma~6.2.
 Given $x=(x_j,\dotsc,x_q)\in\Z^{q-j+1}$,
 we have $L(\eta^j(x))=\rho(b_j)^{-x_j}$.
 Introduce the cultural morphism
 $$
 \theta':
 \U
 \xto{L\Cro\id}
 \U\cro\U
 \xto{M_\U}
 \U
 \xto{\theta^{j+1}}
 \Q^{q-j}_{s_{j+1}\dotso s_q}.
 $$
 Hereafter,
 $\Cro$ combines
 two morphisms with one source
 into a morphism to the product of their targets.
 We have
 $$
 \theta'(\eta^j(x))=
 \theta^{j+1}(L(\eta^j(x))\eta^j(x))=
 \theta^{j+1}(\rho(b_j)^{-x_j}\eta^j(x))=
 $$
 $$
 =
 \theta^{j+1}(\eta^{j+1}(\r(x)))=
 \r(x).
 $$
 Put
 $$
 \theta^j:
 \U
 \xto{\phi_j!\,\Cro\,\theta'}
 \Q_{s_j}\cro\Q^{q-j}_{s_{j+1}\dotso s_q}=
 \Q^{q-j+1}_{s_j\dotso s_q}.
 $$
 We get
 $$
 \theta^j(\eta^j(x))=
 (\phi_j!\,(\eta^j(x)),\theta'(\eta^j(x)))=
 (x_j,\r(x))=
 x.
 \QED
 $$
 \end {demo}

 \begin {claim} [7.2. Lemma.]
 Define a function $m^j$ by the commutative diagram
 $$
 \xymatrix {
 \Z^{q-j+1}_{s_j\dotso s_q}\cro
 \Z^{q-j+1}_{s_j\dotso s_q}
 \ar[rr]^-{m^j}
 \ar[d]_-{\beta^j\cro\beta^j}
 &&
 \Z^{q-j+1}_{s_j\dotso s_q}
 \ar[d]^-{\beta^j}
 \\
 \N^j\cro\N^j
 \ar[rr]^-{M_{\N^j}}
 &&
 \N^j
 .
 }
 $$
 Then
 $m^j$ is a cultural morphism.
 It extends to a cultural morphism
 $$
 \^m^j:
 \Q^{q-j+1}_{s_j\dotso s_q}\cro
 \Q^{q-j+1}_{s_j\dotso s_q}
 \to
 \Q^{q-j+1}_{s_j\dotso s_q}.
 $$
 \end {claim}

 The coordinate polynomials of $m^j$ are known as
 the ``multiplication polynomials''
 \cite[4.2.2]{CMZ}.

 \begin {demo} [Proof.]
 We have the commutative diagram
 $$
 \xymatrix {
 \Z^{q-j+1}_{s_j\dotso s_q}\cro
 \Z^{q-j+1}_{s_j\dotso s_q}
 \ar[rr]^-{m^j}
 \ar[d]_-{\beta^j\cro\beta^j}
 &&
 \Z^{q-j+1}_{s_j\dotso s_q}
 \ar[d]^-{\beta^j}
 \\
 \N^j\cro\N^j
 \ar[rr]^-{M_{\N^j}}
 \ar[d]_-{\rho^j\cro\rho^j}
 &&
 \N^j
 \ar[d]^-{\rho^j}
 \\
 \U\cro\U
 \ar[rr]^-{M_\U}
 &&
 \U
 .
 }
 $$
 Here
 $M_\U$ is a cultural morphism by Lemma~6.1.
 It follows from Lemma~7.1 that
 the composition in the left column is a cultural morphism.
 Since
 the composition in the right column is a cultural immersion
 by Lemma~7.1,
 $m^j$ is a cultural morphism.

 The extension $\^m^j$ exists automatically.
 \qed
 \end {demo}

 \begin {claim} [7.3. Lemma.]
 Given an element $h\in\N^j$,
 define a function $e^j_h$ by the commutative diagram
 $$
 \xymatrix {
 &&
 \Z^{q-j+1}_{s_j\dotso s_q}
 \ar[d]^-{\beta^j}
 \\
 \Z_{s_j}
 \ar[rr]^-{x\mapsto h^x}
 \ar[urr]^-{e^j_h}
 &&
 \N^j
 .
 }
 $$
 Then
 $e^j_h$ is a cultural morphism.
 It extends to a cultural morphism
 $$
 \^e^j_h:
 \Q_{s_j}
 \to
 \Q^{q-j+1}_{s_j\dotso s_q}.
 $$
 \end {claim}

 The coordinate polynomials of $e^j_h$ are a specialization of
 the ``exponentiation polynomials''
 \cite[4.2.2]{CMZ}.

 \begin {demo} [Proof.]
 The composition
 $$
 \Z_{s_j}
 \xto{e^j_h}
 \Z^{q-j+1}_{s_j\dotso s_q}
 \xto{\beta^j}
 \N^j
 \xto{\rho^j}
 \U
 $$
 sends $x$ to $\rho(h)^x$
 and
 thus
 coincides with $E_{\rho(h)}$,
 which is a cultural morphism
 by Lemma~6.2.
 Since $\rho^j\circ\beta^j$ here is a cultural immersion
 by Lemma~7.1,
 $e^j_h$ is a cultural morphism.

 The extension $\^e^j_h$ exists automatically.
 \qed
 \end {demo}


 \head {\S~8. Cultural view of a subgroup of $\N$}

 Let $K\subseteq\N$ be a subgroup.
 For each $j\in\1q$,
 the image of $\N^j\cap K$ in the quotient $\N^j/\N^{j+1}$
 is generated by $\-b_j^{d_j}$
 for some $d_j\ge0$.

 \begin {claim} [8.1. Lemma.]
 There exists
 a cultural morphism $f:\Q^q_{s_1\dotso s_q}\to\Q^q_{s_1\dotso s_q}$ such that
 $$
 \beta^{-1}(K)=
 f^{-1}(d_1\Z\cro\dotso\cro d_q\Z)
 $$
 as subsets of $\Q^q$.
 \end {claim}

 The morhism $f$
 constructed in the proof
 is a cultural isomorphism,
 its $j$th coordinate $f_j$ depends on the first $j$ coordinates of the argument only,
 and,
 for $x\in\Z^q$,
 $f_j(x)\in\Z$ if
 $f_k(x)\in d_k\Z$
 for all $k<j$.
 The proof of these properties is omitted.

 \begin {demo} [Proof.]
 For each $j\in\1{q+1}$,
 we construct a cultural morphism $f^j:\Q^{q-j+1}_{s_j\dotso s_q}\to\Q^{q-j+1}_{s_j\dotso s_q}$ such that
 $$
 (\beta^j)^{-1}(\N^j\cap K)=
 (f^j)^{-1}(d_j\Z\cro\dotso\cro d_q\Z)
 $$
 as subsets of $\Q^{q-j+1}$.
 Backward induction on $j$.
 Let $f^{q+1}:\Q^0\to\Q^0$ be the unique function.
 Take $j\le q$.

 {\it Case $d_j=0$.}
 Then
 $\N^j\cap K\subseteq\N^{j+1}$.
 Put
 $$
 f^j:
 \Q^{q-j+1}_{s_j\dotso s_q}=
 \Q_{s_j}\cro\Q^{q-j}_{s_{j+1}\dotso s_q}
 \xto{\id\cro f^{j+1}}
 \Q_{s_j}\cro\Q^{q-j}_{s_{j+1}\dotso s_q}=
 \Q^{q-j+1}_{s_j\dotso s_q}.
 $$
 Take $x=(x_j,\dotsc,x_q)\in\Q^{q-j+1}$.
 We have
 $$
 (x\in\Z^{q-j+1},
 \
 \beta^j(x)\in\N^j\cap K)
 \quad\Leftrightarrow
 $$
 $$
 \Leftrightarrow\quad
 (x\in0\cro\Z^{q-j},
 \
 \beta^{j+1}(\r(x))\in\N^{j+1}\cap K)
 \quad\Leftrightarrow
 $$
 $$
 \Leftrightarrow\quad
 (x_j=0,
 \
 f^{j+1}(\r(x))\in d_{j+1}\Z\cro\dotso\cro d_q\Z)
 \quad\Leftrightarrow
 $$
 $$
 \quad\Leftrightarrow\quad
 f^j(x)\in0\cro d_{j+1}\Z\cro\dotso\cro d_q\Z.
 $$

 {\it Case $d_j\ne0$.}
 Choose an element $k\in\N^j\cap K$ such that
 $\-k=\-b_j^{d_j}$
 in $\N^j/\N^{j+1}$.
 Consider the cultural morphisms
 $$
 l:
 \Q^{q-j+1}_{s_j\dotso s_q}
 \xto{-\p/d_j}
 \Q_{s_j}
 \xto{\^e^j_k}
 \Q^{q-j+1}_{s_j\dotso s_q},
 $$
 where
 $-\p/d_j:(x_j,\dotsc,x_q)\mapsto-x_j/d_j$
 and
 $\^e^j_k$ is given by Lemma~7.3,
 and
 $$
 f':
 \Q^{q-j+1}_{s_j\dotso s_q}
 \xto{l\Cro\id}
 \Q^{q-j+1}_{s_j\dotso s_q}\cro
 \Q^{q-j+1}_{s_j\dotso s_q}
 \xto{\^m^j}
 \Q^{q-j+1}_{s_j\dotso s_q}
 \xto\r
 \Q^{q-j}_{s_{j+1}\dotso s_q}
 \xto{f^{j+1}}
 \Q^{q-j}_{s_{j+1}\dotso s_q},
 $$
 where
 $\^m^j$ is given by Lemma~7.2.
 Put
 $$
 f^j:
 \Q^{q-j+1}_{s_j\dotso s_q}
 \xto{\p\Cro f'}
 \Q_{s_j}\cro\Q^{q-j}_{s_{j+1}\dotso s_q}=
 \Q^{q-j+1}_{s_j\dotso s_q}.
 $$

 Take $x=(x_j,\dotsc,x_q)\in d_j\Z\cro\Z^{q-j}$.
 Then
 $$
 k^{-x_j/d_j}\beta^j(x)\in\N^{j+1}
 $$
 and
 $$
 k^{-x_j/d_j}\beta^j(x)=\beta^j(y),
 $$
 where
 $y\in\Z^{q-j+1}$,
 $$
 y=
 \^m^j(\^e^j_k(-x_j/d_j),x).
 $$
 Thus
 $\p(y)=0$
 and
 $$
 k^{-x_j/d_j}\beta^j(x)=
 \beta^j(y)=
 \beta^{j+1}(\r(y))=
 \beta^{j+1}(\r(\^m^j(\^e^j_k(-x_j/d_j),x))).
 $$

 Take $x=(x_j,\dotsc,x_q)\in\Q^{q-j+1}$.
 Put
 $y=\^m^j(\^e^j_k(-x_j/d_j),x)\in\Q^{q-j+1}$
 and
 $y'=\r(y)\in\Q^{q-j}$.
 We show that
 \begin {equation} \label {py}
 \p(y)=0
 \end {equation}
 and
 \begin {equation} \label {x}
 x=
 \^m^j(\^e^j_k(x_j/d_j),y).
 \end {equation}
 It suffices to consider the case $x\in d_j\Z\cro\Z^{q-j}$.
 Then,
 as shown above,
 $y\in\Z^{q-j+1}$,
 $\p(y)=0$,
 and
 $$
 k^{-x_j/d_j}\beta^j(x)=\beta^j(y).
 $$
 Thus
 $$
 \beta^j(x)=k^{x_j/d_j}\beta^j(y),
 $$
 which implies \eqref{x}.
 It follows from
 \eqref{py}
 and
 \eqref{x}
 that
 $$
 (x_j\in d_j\Z,
 \
 y'\in\Z^{q-j})
 \quad\Rightarrow\quad
 x\in d_j\Z\cro\Z^{q-j}.
 $$
 We have
 $f'(x)=f^{j+1}(y')$
 and
 $$
 (x\in\Z^{q-j+1},
 \
 \beta^j(x)\in\N^j\cap K)
 \quad\Leftrightarrow
 $$
 $$
 \Leftrightarrow\quad
 (x\in d_j\Z\cro\Z^{q-j},
 \
 k^{-x_j/d_j}\beta^j(x)\in\N^{j+1}\cap K)
 \quad\Leftrightarrow
 $$
 $$
 \Leftrightarrow\quad
 (x_j\in d_j\Z,
 \
 y'\in\Z^{q-j},
 \
 \beta^{j+1}(y')\in\N^{j+1}\cap K)
 \quad\Leftrightarrow
 $$
 $$
 \Leftrightarrow\quad
 (x_j\in d_j\Z,
 \
 f^{j+1}(y')\in d_{j+1}\Z\cro\dots\cro d_q\Z)
 \quad\Leftrightarrow
 $$
 $$
 \Leftrightarrow\quad
 (x_j\in d_j\Z,
 \
 f'(x)\in d_{j+1}\Z\cro\dotso\cro d_q\Z)
 \quad\Leftrightarrow
 $$
 $$
 \Leftrightarrow\quad
 f^j(x)\in d_j\Z\cro\dotso\cro d_q\Z.
 \QED
 $$
 \end {demo}


 \head {\S~9. Defining $\gamma^{r+1}G$ by equations and congruences}

 \begin {claim} [9.1. Theorem.]
 Let
 $G$ be a group with elements $g_1,\dotsc,g_n\in G$.
 Fix $r\ge0$.
 Then,
 for some $q\ge0$,
 there are
 rational polynomials $P_j(X_1,\dotsc,X_n)$,
 $j\in\1q$,
 of degree at most $r$
 and
 integers $d_1,\dotsc,d_q\ge0$
 such that,
 for any $x=(x_1,\dotsc,x_n)\in\Z^n$,
 $$
 g_1^{x_1}\dotso g_n^{x_n}\in\gamma^{r+1}G
 \quad\Leftrightarrow\quad
 (P_j(x)\in d_j\Z,
 \
 j\in\1q).
 $$
 \end {claim}

 The polynomials $P_j$
 constructed in the proof
 have the following integrality property:
 for $x\in\Z^n$,
 $P_j(x)\in\Z$ if
 $P_k(x)\in d_k\Z$ for $k<j$.
 The check is omitted.
 
 \begin {demo} [Proof.]
 We use the constructions of the previous sections of Part~II
 for the given $n$ and $r$.
 In particular,
 we let the required $q$ be the size of the Mal\6cev basis of $\N$,
 the free nilpotent group of rank $n$ and class $r$.
 Consider the homomorphism
 $$
 t:\N\to G/\gamma^{r+1}G,
 \qquad
 \-Z_i\mapsto\-g_i.
 $$
 Put $K=\Ker t\subseteq\N$.
 By Lemma~8.1,
 there are
 integers $d_1,\dotsc,d_q\ge0$
 and
 a cultural morphism $f:\Q^q_{s_1\dotso s_q}\to\Q^q_{s_1\dotso s_q}$
 such that
 $$
 \beta^{-1}(K)=f^{-1}(d_1\Z\cro\dotso\cro d_q\Z)
 $$
 as subsets of $\Q^q$.
 Define the required polynomials $P_j$ by the equality
 $$
 f(x_1,\dotsc,x_n,0,\dotsc,0)=
 (P_j(x))_{j\in\1q},
 \qquad
 x=(x_1,\dotsc,x_n)\in\Q^n.
 $$
 Since $s_1=\dotso=s_n=1$,
 the degree of $P_j$ is at most $s_j$,
 which does not exceed $r$.
 Given $x=(x_1,\dotsc,x_n)\in\Z^n$,
 we have
 $$
 g_1^{x_1}\dotso g_n^{x_n}\bmod\gamma^{r+1}G=
 t(\-Z_1^{x_1}\dotso\-Z_n^{x_n})=
 t(\beta(x_1,\dotsc,x_n,0,\dotsc,0))
 $$
 in $G/\gamma^{r+1}G$
 and thus
 $$
 g_1^{x_1}\dotso g_n^{x_n}\in\gamma^{r+1}G
 \quad\Leftrightarrow\quad
 \beta(x_1,\dotsc,x_n,0,\dotsc,0)\in K
 \quad\Leftrightarrow
 $$
 $$
 \Leftrightarrow\quad
 f(x_1,\dotsc,x_n,0,\dotsc,0)\in d_1\Z\cro\dotso\cro d_q\Z
 \quad\Leftrightarrow\quad
 (P_j(x)\in d_j\Z,
 \
 j\in\1q).
 \QED
 $$
 \end {demo}


 \suphead {Part III}


 In this part,
 we consider the group $[S^1,Y]=\pi_1(Y)$.


 \head {\S~10. Managing an ensemble of maps $S^1\to Y$}

 For
 $n\ge0$ and
 a group $G$,
 introduce the function
 $$
 M:
 G^n\to G,
 \qquad
 (g_1,\dotsc,g_n)\mapsto g_1\dotso g_n.
 $$
 For $K\subseteq\1n$,
 let $\omega_K:G^n\to G^K$ be the projection.

 \begin {claim} [10.1. Lemma.]
 Consider an ensemble $A\in\<Y^{S^1}\>$,
 $$
 A=
 \sum_{i\in I}
 u_i\`a_i\',
 $$
 such that
 $A\%=r0$.
 Then,
 for some $n\ge1$,
 there exist elements $z_i\in\pi_1(Y)^n$, $i\in I$, such that
 $[a_i]=M(z_i)$ and
 the element
 \begin {equation} \label {Z}
 Z=
 \sum_{i\in I}
 u_i\`z_i\'
 \in\<\pi_1(Y)^n\>
 \end {equation}
 satisfies $\<\omega_K\>(Z)=0$ in $\<\pi_1(Y)^K\>$
 for all $K\subseteq\1n$
 with $|K|\le r$.
 \end {claim}

 \begin {demo} [Proof.]
 Take a finite subspace $D\subseteq S^1$
 consisting of $n\ge2$ points.
 It cuts $S^1$ into closed arcs $B_k$, $k\in\1n$.
 A continuous function $v:B_k\to Y$
 with $v(\partial B_k)=\{\0_Y\}$
 has the
 (relative to $\partial B_k$)
 homotopy class $[v]\in\pi_1(Y)$.
 For a map $w:S^1\to Y$
 with $w(D)=\{\0_Y\}$,
 we have
 \begin {equation} \label {pro}
 [w]=
 \prod_{k=1}^n
 [w|_{B_k}]
 \end {equation}
 in $\pi_1(Y)$
 (we assume that $B_k$ are oriented and numbered properly).

 By \cite[Corollary~6.2]{sim},
 we may assume that
 $A\$=r\Gamma0$
 for some open cover $\Gamma$ of $S^1$.
 We suppose that
 $D$ is chosen dense enough so that
 each $B_k$ is contained in some $G_k\in\Gamma$.
 Put
 $$
 V=
 \bigbou_{i\in I}
 S^1.
 $$
 Let $U$ be the quotient of $V$ by the identifications
 $\ins_i(x)\approx\ins_j(x)$
 for
 $x\in B_k$ and
 $i,j\in I$ such that
 $a_i\|={G_k}a_j$.
 $U$ is a graph.
 Let $h:V\to U$ be the projection.
 Introduce the maps
 $$
 e_i:
 S^1
 \xto{\ins_i}
 V
 \xto{h}
 U.
 $$
 There is a map $q:U\to S^1$ such that
 $q\circ e_i=\id_{S^1}$.
 Put
 $$
 b=
 \bigBou_{i\in I}
 a_i
 :V\to Y.
 $$
 There is a map $a:U\to Y$ such that
 $b=a\circ h$.
 Clearly,
 $a_i=a\circ e_i$.
 Put $\~D=q^{-1}(D)\subseteq U$.
 $\~D$ is a finite subspace.
 The map $a|_{\~D}$ is null-homotopic because
 the inclusion $\~D\to U$ is.
 Extending the homotopy,
 we get a map $\^a:U\to Y$ such that
 $\^a\sim a$ and
 $\^a(\~D)=\{\0_Y\}$.
 Put $\^a_i=\^a\circ e_i:S^1\to Y$.
 Clearly,
 $\^a_i\sim a_i$ and
 $\^a_i(D)=\{\0_Y\}$.
 Put
 $$
 z_i=
 ([\^a_i|_{B_k}])_{k\in\1n}
 \in\pi_1(Y)^n.
 $$
 We have
 $$
 [a_i]=
 [\^a_i]
 \overset{(*)}=
 \prod_{k=1}^n
 [\^a_i|_{B_k}]=
 M(z_i),
 $$
 where
 ($*$) follows from \eqref{pro}.
 For $k\in\1n$ and $i,j\in I$,
 we have the implication
 \begin {equation} \label {ergo}
 a_i\|={G_k}a_j
 \quad
 \Rightarrow
 \quad
 [\^a_i|_{B_k}]=[\^a_j|_{B_k}]
 \end {equation}
 because the premise implies that
 $e_i\|={B_k}e_j$ and
 thus
 $\^a_i\|={B_k}\^a_j$.

 Consider the element $Z\in\<\pi_1(Y)^n\>$ given by \eqref{Z}.
 Take $K\subseteq\1n$.
 Put
 $$
 G(K)=
 \{\0_Y\}\cup
 \bigcup_{k\in K}
 G_k
 \subseteq S^1.
 $$
 By \eqref{ergo},
 we have the implication
 $$
 a_i\|={G(K)}a_j
 \quad
 \Rightarrow
 \quad
 \omega_K(z_i)=\omega_K(z_j).
 $$
 Suppose that
 $|K|\le r$.
 Then
 $A|_{G(K)}=0$ because
 $A\$=r\Gamma0$.
 Thus
 $\<\omega_K\>(Z)=0$.
 \qed
 \end {demo}


 \head {\S~11. Similarity on $\pi_1(Y)$}

 \begin {claim} [11.1. Lemma.]
 Let $G$ be a group.
 Consider an element $Z\in\<G^n\>$,
 $$
 Z=
 \sum_{i\in I}
 u_i\`z_i\',
 $$
 where
 $I$ has a distinguished element $0$ and
 $u_0=1$.
 Suppose that
 $\<\omega_K\>(Z)=0$
 for all $K\subseteq\1n$
 with $|K|\le r$ and
 $M(z_i)\in\gamma^{r+1}G$
 for all $i\ne0$.
 Then
 $M(z_0)\in\gamma^{r+1}G$.
 \end {claim}

 \begin {demo} [Proof.]
 We have $z_i=(z_{i1},\dotso,z_{in})$,
 where $z_{ik}\in G$.
 Take distinct $g_1,\dotsc,g_m\in G$
 that include all the $z_{ik}$.
 We have
 $$
 z_{ik}=
 \prod_{l=1}^m
 g_l^{\[z_{ik}=g_l\]}
 $$
 and
 thus
 $$
 M(z_i)=
 \prod_{k=1}^n
 \prod_{l=1}^m
 g_l^{\[z_{ik}=g_l\]}.
 $$
 Hereafter,
 given a condition $C$,
 the integer $\[C\]$ is
 $1$ under $C$ and
 $0$ otherwise.
 By Theorem~9.1,
 for some $q\ge0$,
 there are
 rational polynomials $P_j(X)$,
 $X=(X_{kl})_{k\in\1n,\,l\in\1m}$,
 $j\in\1q$,
 of degree at most $r$ and
 integers $d_j\ge0$
 such that,
 for any collection $x=(x_{kl})_{k\in\1n,\,l\in\1m}$, $x_{kl}\in\Z$,
 we have the equivalence
 $$
 \prod_{k=1}^n
 \prod_{l=1}^m
 g_l^{x_{kl}}
 \in\gamma^{r+1}G
 \quad
 \Leftrightarrow 
 \quad
 (P_j(x)\in d_j\Z,\ 
 j\in\1q).
 $$
 Order the set $\1n\times\1m$ totally.
 We have
 $$
 P_j(X)=
 \sum_{\substack{
         0\le s\le r, \\
         (k_1,l_1)\le\dotso\le(k_s,l_s)
 }}
 P^{(s)}_{jk_1l_1\dotso k_sl_s}
 X_{k_1l_1}\dotso X_{k_sl_s}
 $$
 for some $P^{(s)}_{jk_1l_1\dotso k_sl_s}\in\Q$.
 We have
 \begin {multline} \label {summa}
 \sum_{i\in I}
 u_iP_j((\[z_{ik}=g_l\])_{k\in\1n,\,l\in\1m})= \\ =
 \sum_{\substack{
         0\le s\le r, \\
         (k_1,l_1)\le\dotso\le(k_s,l_s)
 }}
 P^{(s)}_{jk_1l_1\dotso k_sl_s}
 \sum_{i\in I}
 u_i\[z_{ik_t}=g_{l_t},\ t\in\1s\]
 \overset{(*)}=
 0,
 \end {multline}
 where
 ($*$) holds because
 the inner sum is zero,
 which is because
 $\<\omega_K\>(Z)=0$
 for $K=\{k_1,\dotsc,k_s\}$.
 Since $M(z_i)\in\gamma^{r+1}(G)$
 for $i\ne0$,
 we have
 \begin {equation} \label {cong}
 P_j((\[z_{ik}=g_l\])_{k\in\1n,\,l\in\1m})
 \in d_j\Z
 \end {equation}
 for $i\ne0$.
 Since $u_0=1$,
 it follows from \eqref{summa} that
 \eqref{cong} holds for $i=0$ too.
 Thus
 $M(z_0)\in\gamma^{r+1}(G)$.
 \qed
 \end {demo}

 \begin {claim} [11.2. Theorem.]
 Let $Y$ be a cellular space.
 Then
 \begin {equation} \label {main}
 \pi_1(Y)^{(r+1)}=\gamma^{r+1}\pi_1(Y).
 \end {equation}
 \end {claim}

 \begin {demo} [Proof.]
 The inclusion $\supseteq$ in \eqref{main}
 follows from Theorem~4.3.
 To prove the inclusion $\subseteq$,
 we take $a\in Y^{S^1}$ such that
 $\0\%\sim r a$
 and
 check that
 $[a]\in\gamma^{r+1}\pi_1(Y)$.
 There is an ensemble $D\in\<Y^{S^1}\>$,
 $$
 D=
 \sum_i
 u_u\`d_i\',
 $$
 where $d_i\sim\0$,
 such that
 $D\%=r\`a\'$.
 By Lemma~10.1,
 for some $n\ge1$,
 there are elements $z,w_i\in\pi_1(Y)^n$ such that
 $M(z)=[a]$ and
 $M(w_i)=1$
 in $\pi_1(Y)$ and,
 putting
 $$
 W=
 \sum_i
 u_i\`w_i\'
 \in\<\pi_1(Y)^n\>,
 $$
 we have $\<\omega_K\>(\`z\'-W)=0$
 for all $K\subseteq\1n$
 with $|K|\le r$.
 By Lemma~11.1,
 $M(z)\in\gamma^{r+1}\pi_1(Y)$,
 which is what we need.
 \qed
 \end {demo}


 \head {\S~12. Finite-order invariants on $\pi_1(Y)$}

 For a group $G$,
 $\<G\>$ is its group ring.
 Let $[G]\subseteq\<G\>$ be the augmentation ideal,
 i.~e.,
 the kernel of the ring homomorphism
 (called the augmentation)
 $\<G\>\to\Z$,
 $\`g\'\mapsto1$.

 \begin {claim} [12.1. Lemma.]
 Let
 $G$ be a group and
 $Z\in\<G^n\>$ be an element such that
 $\<\omega_K\>(Z)=0$
 in $\<G^K\>$
 for all $K\subseteq\1n$
 with $|K|\le r$.
 Then
 $\<M\>(Z)\in[G]^{r+1}$
 ($\subseteq\<G\>$).
 \end {claim}

 \begin {demo} [Proof.]
 For $K\subseteq\1n$,
 consider the function
 $$
 \epsilon_K:G^K\to G^n,
 \qquad
 (g_k)_{k\in K}
 \mapsto
 (\~g_k)_{k\in\1n},
 $$
 where
 $\~g_k$ equals
 $g_k$ if $k\in K$
 and
 $1$ otherwise,
 the composition
 $$
 \rho_K:
 G^n
 \xto{\omega_K}
 G^K
 \xto{\epsilon_K}
 G^n
 $$
 and
 the homomorphism $S_K:\<G^n\>\to\<G^n\>$,
 $$
 S_K=
 \sum_{L\subseteq K}
 (-1)^{|L|}\<\rho_L\>.
 $$
 If $K=\{k_1,\dotsc,k_t\}$,
 $k_1<\dotso<k_t$,
 then
 $$
 (\<M\>\circ S_K)
 (\`(g_k)_{k\in\1n}\')=
 (1-\`g_{k_1}\')\dotso(1-\`g_{k_t}\')
 $$
 in $\<G\>$.
 Thus
 \begin {equation} \label {img}
 \Im(\<M\>\circ S_K)\subseteq[G]^{|K|}.
 \end {equation}
 We have
 \begin{multline*}
 \sum_{K\subseteq\1n}
 (-1)^{|K|}S_K=
 \sum_{K\subseteq\1n}
 (-1)^{|K|}
 \sum_{L\subseteq K}
 (-1)^{|L|}\<\rho_L\>=
 \\
 =
 \sum_{L\subseteq\1n}
 (-1)^{|L|}
 \bigl(
 \sum_{K\subseteq\1n:K\supseteq L}
 (-1)^{|K|}
 \bigr)
 \<\rho_L\>.
 \end{multline*}
 The inner sum equals $(-1)^n\[L=\1n\]$.
 Thus
 $$
 \sum_{K\subseteq\1n}
 (-1)^{|K|}S_K=
 \<\rho_{\1n}\>=
 \id_{\<G\>}.
 $$
 For $L\subseteq\1n$,
 $|L|\le r$,
 we have
 $
 \<\rho_L\>(Z)=
 \<\epsilon_L\>(\<\omega_L\>(Z))=
 0
 $.
 Thus
 $S_K(Z)=0$
 if $|K|\le r$.
 We get
 $$
 Z=
 \sum_{K\subseteq\1n}
 (-1)^{|K|}S_K(Z)=
 \sum_{K\subseteq\1n:|K|\ge r+1}
 (-1)^{|K|}S_K(Z).
 $$
 Thus
 $$
 \<M\>(Z)=
 \sum_{K\subseteq\1n:|K|\ge r+1}
 (-1)^{|K|}(\<M\>\circ S_K)(Z).
 $$
 By \eqref{img},
 $\<M\>(Z)\in[G]^{r+1}$.
 \qed
 \end {demo}

 A function $f:G\to L$,
 where $L$ is an abelian group,
 gives rise to the homomorphism
 $$
 \pp f:\<G\>\to L,
 \qquad
 \`g\'\mapsto f(g).
 $$
 We define $\deg f\in\{-\infty,0,1,\dotsc,\infty\}$,
 the {\it degree\/} of $f$,
 as the infimum of $r\in\Z$ such that
 $\pp f|_{[G]^{r+1}}=0$
 (adopting $[G]^s=\<G\>$ for $s\le0$).

 \begin {claim} [12.2. Theorem.]
 Let
 $Y$ be a cellular space,
 $L$ be an abelian group and
 $f:\pi_1(Y)\to L$
 be a homotopy invariant
 (i.~e., a function).
 Then
 $\ord f=\deg f$
 \end {claim}

 \begin {demo} [Proof.]
 We suppose $f\ne0$
 omitting the converse case.

 (1)
 Suppose that
 $\ord f\le r$
 ($r\ge0$).
 We show that
 $\deg f\le r$.
 It suffices to check that
 $$
 \pp f((1-\`[a_1]\')\dotso(1-\`[a_{r+1}]\'))=0
 $$
 for any $a_1,\dotsc,a_{r+1}\in Y^{S^1}$.
 Put
 $W=S^1\bou\dotso\bou S^1$
 ($r+1$ summands)
 and
 $$
 q=
 a_1\Bou\dotso\Bou a_{r+1}:
 W\to Y.
 $$
 Let
 $p:S^1\to W$ be the $(r+1)$-fold comultiplication
 and
 $\Lambda_d:W\to W$,
 $d\in\E^{r+1}$,
 be as in \cite[\S~3]{sim}.
 Consider the ensemble $A\in\<Y^{S^1}\>$,
 $$
 A=
 \sum_{d\in\E^{r+1}}
 (-1)^{|d|}
 \`a(d)\',
 $$
 where
 $$
 a(d):
 S^1
 \xto{p}
 W
 \xto{\Lambda(d)}
 W
 \xto{q}
 Y.
 $$
 Clearly,
 $$
 [a(d)]=
 [a_1]^{d_1}\dotso[a_{r+1}]^{d_{r+1}}
 $$
 in $\pi_1(Y)$.
 By \cite[Lemma~3.1]{sim},
 $A\%=r0$.
 We have
 \begin {multline*}
 \pp f((1-\`[a_1]\')\dotso(1-\`[a_{r+1}]\'))=
 \sum_{d\in\E^{r+1}}
 (-1)^{|d|}
 f([a_1]^{d_1}\dotso[a_{r+1}]^{d_{r+1}})= \\ =
 \sum_{d\in\E^{r+1}}
 (-1)^{|d|}
 f([a(d)])
 \overset{(*)}=
 0,
 \end {multline*}
 where
 ($*$) holds because
 $\ord f\le r$.

 (2)
 Suppose that
 $\deg f\le r$
 ($r\ge0$).
 We show that
 $\ord f\le r$.
 Take an ensemble $A\in\<Y^{S^1}\>$,
 $$
 A=
 \sum_{i\in I}
 u_i\`a_i\',
 $$
 such that
 $A\%=r0$.
 We should show that
 $$
 \sum_{i\in I}
 u_if([a_i])=
 0.
 $$
 By Lemma~10.1,
 for some $n\ge1$,
 there exist elements $z_i\in\pi_1(Y)^n$, $i\in I$, such that
 $[a_i]=M(z_i)$ and
 the element $Z\in\<\pi_1(Y)^n\>$ given by \eqref{Z}
 satisfies $\<\omega_K\>(Z)=0$ in $\<\pi_1(Y)^K\>$
 for all $K\subseteq\1n$
 with $|K|\le r$.
 We have
 $$
 \sum_{i\in I}
 u_if([a_i])=
 \pp f(\<M\>(Z)).
 $$
 By Lemma~12.1,
 $\<M\>(Z)\in[G]^{r+1}$.
 Since $\deg f\le r$,
 $\pp f(\<M\>(Z))=0$.
 \qed
 \end {demo}



 \noindent
 \href{mailto:ssp@pdmi.ras.ru}{\tt ssp@pdmi.ras.ru}

 \noindent
 \url{http://www.pdmi.ras.ru/~ssp}

 \end {document}